\documentclass{amsart}
\usepackage{amsmath,amsthm,amssymb,IMjournal}
\usepackage{graphicx}

\newcommand{\bg}{\mathbf{g}}
\newcommand{\bx}{\mathbf{x}}
\newcommand{\bjedna}{\mathbf{1}}

\newcommand{\bu}{\mathbf{u}}

\newcommand{\Zosem}{\mathbb{Z}_8}

\DeclareMathOperator{\rng}{rng}

\DeclareMathOperator{\od}{mod}
\DeclareMathOperator{\Sub}{Sub}

\begin{document}
\newtheorem{theorem}{Theorem}[section]
\newtheorem{lemma}{Lemma}[section]
\newtheorem{definition}[lemma]{Definition}
\newtheorem{corollary}[lemma]{Corollary}

\numberwithin{equation}{section}

\title{EXPANSIONS OF THE GROUP $\Zosem$ (PART I)}

\author{\|Miroslav |Ploščica|,
\|Radka |Schwartzová|, 
\|Ivana |Varga|}

\rec {\today}


\abstract 
We investigate clones in the interval between the group polynomials and the ring polynomials of $\Zosem$. This is the simplest open case of the problem, as the answer is known for
${\mathbb Z}_{p^2}$ (with $p$ prime) and, in general, ${\mathbb Z}_n$ reduces to the case when $n$ is a prime power. The investigated structure proves to be very complicated,
so we provide only a partial description. We restrict our attention to polynomials
whose nonlinear monomials have even coefficients.
\endabstract



\thanks
   {}The research has been supported by VEGA Grant 1/0152/22 and the EU NextGenerationEU through the Recovery and Resilience Plan for Slovakia under the project No. 09I03-03-V05-00008.
\endthanks

\section{Introduction}\label{sec1}

 A family of finitary operations on a set $A$ that contains all projections and is closed under composition is called a clone.
 Clones have an important position in~modern universal algebra. This is due to the fact that many important properties of algebraic structures are determined by their clones of term or polynomial operations.

All clones on a given set $A$ form a complete lattice.
For a $2$-element set $A$ this lattice is countable and has been described by E. Post \cite{post}. 
If the considered set has at least three elements, the lattice of clones is uncountable and there is no chance to describe it completely, see \cite{goldstern2008survey}. 
Therefore, the research usually aims to~describe some interesting parts of this lattice. Very often, the investigated clones are connected with some special algebraic structure on the base set. For~instance, various papers study clones containing a group operation (\cite{Bulatovclonynagrupe, Bulatovclonynagrupe2, Szendreiclonynagrupe, Salomaaclonynagrupe}).
 The~papers   \cite{Gavrilov1996clonynaokruhu, Gavrilov1997clonynaokruhu, Remizov1989clonynaokruhu, Rosenberg1995clonynaokruhu, rosenberg, bulatclonynaokruhu, AichigerMayr, Mayr} are focused on the study of clones of the ring polynomials of $\mathbb{Z}_n$.  

In our paper we investigate the interval 
$\mathcal{J}_n= \langle P(\mathbb{Z}_n, +), P(\mathbb{Z}_n, +, \cdot) \rangle$ in the~lattice of clones. Here, 
$ P(\mathbb{Z}_n, +)$ is the clone of all polynomial operations of the group $(\mathbb{Z}_n, +)$, while
$P(\mathbb{Z}_n, +, \cdot)$ denotes the clone of all polynomial operations of the ring $(\mathbb{Z}_n, +, \cdot)$.
Hence, the elements of $ P(\mathbb{Z}_n, +)$ are of the  form 
$$p(x_1, \dots, x_k)=a_0+a_1x_1 +a_2x_2 +\dots + a_nx_k,$$
where $a_0 \in \mathbb{Z}_n$, $a_1, \dots ,a_n\in \mathbb{Z}$. 
Operations in $P(\mathbb{Z}_n, +, \cdot)$ can be expressed by polynomials of the form
\begin{equation}\label{naryf}f(\bx)=\sum_\alpha a_\alpha\bx^\alpha,\end{equation}
where $\bx=(x_1,\dots, x_n)$, the sum is taken over finitely many $n$-tuples $\alpha=(\alpha_1,\dots,\alpha_n)$ of natural numbers, the coefficients $a_\alpha$
belong to $\mathbb{Z}_n$, and $\bx^\alpha$ is an abbreviation for $x_1^{\alpha_1}\dots x_n^{\alpha_n}$.

\begin{lemma}\label{LJPk}
Let $n, m$ be natural numbers. 
If $(m,n)=1$, then 
$$\mathcal{J}_{mn} \cong \mathcal{J}_m \times \mathcal{J}_n.$$

\end{lemma}
The proof of the Lemma \ref{LJPk} is implicitly included in  \cite{ploscica2020clones}. For the sake of completeness, we provide the proof in the next section.
This lemma means that in order to describe the interval $\mathcal{J}_{n}$, it suffices to investigate the case $n=p^k$, where $p$ is a prime number.

The case $k=1$ is well-known. Since $\mathbb{Z}_{p}$ is a field, every operation is a ring polynomial. By I. Rosenberg (\cite{rosenberg}), clone $P(\mathbb{Z}_{p}, +)$ is maximal, which means that the interval $\mathcal{J}_{p}$ consists  of just two clones: $P(\mathbb{Z}_{p}, +)$, and $P(\mathbb{Z}_{p}, +, \cdot)$. 

The case $k=2$ is much more complicated.
The detailed description of $\mathcal{J}_{p^2}$ was provided by A. Bulatov in \cite{bulatclonynaokruhu}.  The lattice is depicted below.

\begin{figure}[h]       
\centering
\includegraphics[width= 0.4\textwidth]{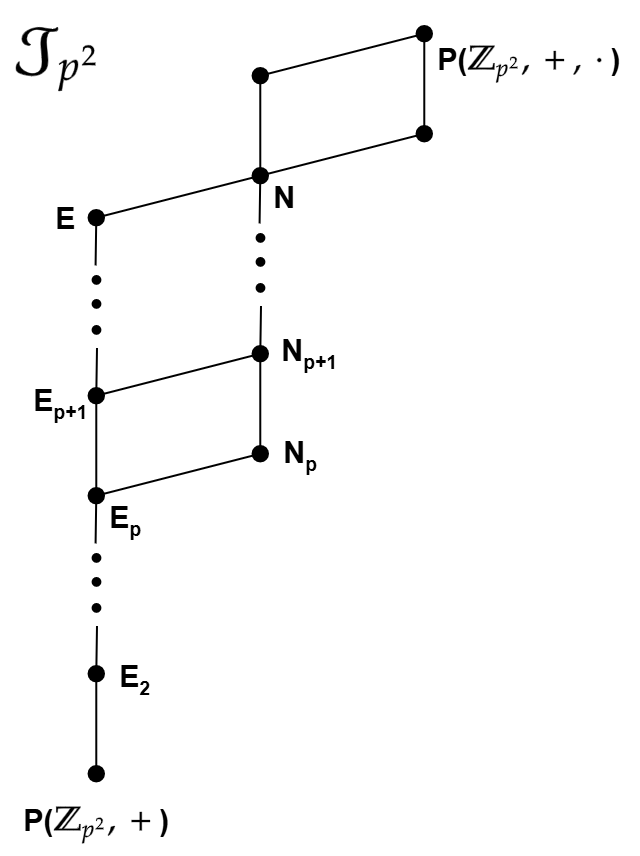}
\caption{The Lattice of Polynomial Clones of $\mathcal{J}_{p^2}$}\label{obr-03}
\end{figure}

(We would like to remark that
$\mathcal{J}_{p^2}$ for $p=2$ was described by Krokhin et. al. \cite{Krokhin1997} and generalized by A. Bulatov and P. Idziak \cite{bulatov2003counting}.)

The case $k=3$ is still open. In our paper, we consider the case $p=2$.
Even this special case is so complicated, that we are able
to describe only part of the interval $\mathcal{J}_{8}$,
namely the interval $\langle P(\mathbb{Z}_{8}, +), M_1 \rangle$, where the clone $M_1$ consists of all ring polynomials with the property that all nonlinear monomials have even coefficients.

\section{Reduction to prime power}

In this section, we prove Lemma \ref{LJPk}. In fact, we show slightly more. 
We regard  ${\mathbb Z}_n$ as the set with elements
$0,\dots,n-1$. Recall that the group $\mathbb{Z}_n$ and
the ring $\mathbb{Z}_n$ have the same congruences, namely
$\od m$ for every $m$ dividing $n$ defined by
$$x\equiv y\ (\od m) \quad\text{iff}\quad m|(x-y).$$
We also use "$mod$" as an operational symbol, so that
$x\mod m$ is the unique element of ${\mathbb Z}_m$ congruent with $x$ modulo $m$.
Let $m$ and $n$ be relatively prime natural numbers. Let $\mathcal{O}$ be the set of all finitary operations on
the set ${\mathbb Z}_{mn}$ that preserve congruences $\od m$ and $\od n$. Further, let $\mathcal{O}_m$
and $\mathcal{O}_n$ be the sets of all finitary operations on ${\mathbb Z}_m$ and ${\mathbb Z}_n$,
respectively.

For any $k$-ary $f\in\mathcal{O}$ we define $f_m\in\mathcal{O}_m$ and $f_n\in\mathcal{O}_n$ by
$$f_m(x_1,\dots,x_k)=f(x_1,\dots,x_k) \quad \od m,$$
$$f_n(x_1,\dots,x_k)=f(x_1,\dots,x_k)\quad\od n.$$
Further, for $k$-ary $(g,h)\in\mathcal{O}_m\times\mathcal{O}_n$ 
and
$x_1,\dots,x_k\in{\mathbb Z}_{mn}$, we define $(g\times h)(x_1,\dots,x_k)$ as the unique element of ${\mathbb Z}_{mn}$
satisfying
$$(g\times h)(x_1,\dots,x_k)\equiv g(x_1\od m,\dots,x_k\od m)\quad (\od m),$$
$$ (g\times h)(x_1,\dots,x_k)\equiv h(x_1\od n,\dots,x_k\od n)\quad (\od n),$$
which exists by the Chinese remainder theorem. Clearly, this operation preserves congruences mod $m$
and mod $n$, so it belongs to ${\mathcal O}$. 

\begin{lemma} \label{iso0} 
For every $k$-ary $f\in{\mathcal O}$, $g\in{\mathcal O}_m$, $h\in{\mathcal O}_n$, the
equalities $(g\times h)_m=g$, $(g\times h)_n=h$ 
and $f_m\times f_n=f$ hold.
\end{lemma}
\proof   The first two equalities are easy. To prove the third one, we need to show that for any $x_1,\dots,x_k\in{\mathbb Z}_{mn}$,
$$ f(x_1,\dots,x_k)\equiv f_m(x_1\od m,\dots,x_k\od m)\quad (\od m),$$
$$ f(x_1,\dots,x_k)\equiv f_n(x_1\od n,\dots,x_k\od n)\quad (\od n).$$
By the definition of $f_m$,
$$f_m(x_1\od m,\dots,x_k\od m)\equiv f(x_1\od m,\dots,x_k\od m)\quad (\od m).$$
Since $f$ preserves the congruence mod $m$, we have 
$$f(x_1\od m,\dots,x_k\od m)\equiv f(x_1,\dots,x_k)\ (\od m),$$
so we obtain the statement for $m$. The argument for $n$ is the same.
\endproof

\begin{lemma} The equality
$$f(g^1,\dots,g^k)_m=f_m(g^1_m,\dots,g^k_m)$$ holds
for every $k$-ary $f$ and $l$-ary $g^i$ in $\mathcal{O}$.\label{compos}
\end{lemma}

\proof Let $k$-ary $f$ and $l$-ary $g^1,\dots,g^k$
belong to ${\mathcal O}$. Let $\bx=(x_1,\dots,x_l)\in ({\mathbb Z}_m)^l$. For every $i$ we have
$$g^i_m(\bx)\equiv g^i(\bx)\quad (\od m).$$
Since $f$ preserves the congruence mod $m$, we have
$$f(g^1(\bx),\dots,g^k(\bx))\equiv f(g^1_m(\bx),\dots,g^k_m(\bx))\quad (\od m).$$
Also,
$$f_m(g^1_m(\bx),\dots,g^k_m(\bx))\equiv f(g^1_m(\bx),\dots,g^k_m(\bx))\quad (\od m),$$
which implies $f(g^1,\dots,g^k)_m=f_m(g^1_m,\dots,g^k_m)$.
\endproof

The correspondence described in Lemma \ref{iso0} extends to subclones as follows. Let $C$ be a subclone of ${\mathcal O}$. We define
$$C_m=\{f_m\mid f\in {\mathcal O}\};$$
$$C_n=\{f_n\mid f\in {\mathcal O}\}.$$

Clearly, if $f$ is a projection, then both $f_m$ and $f_n$ are projections on the same variable. Now Lemma \ref{compos} has the following consequence.

\begin{lemma} For every clone $C\subseteq{\mathcal O}$,  $C_m$ and $C_n$ are subclones of ${\mathcal O}_m$ and ${\mathcal O}_n$, respectively. Moreover, if $C$ is generated by its subset $F$, then $C_m$ and $C_n$
are generated by $\{f_m\mid f\in F\}$ and $\{f_n\mid f\in F\}$, respectively.
\end{lemma}

Further, Let  $D$ and $E$ be subclones of ${\mathcal O}_m$ and ${\mathcal O}_n$, respectively. 
We define
$$D\times E=\{g\times h\mid g\in{\mathcal O}_m,\ h\in{\mathcal O}_n,\ \text{of the same arity}\}.$$

\begin{lemma} For any $D$ and $E$, $D\times E$ is a clone.\end{lemma}
\proof If $g$ and $h$
are projections on the same variable, then $g\times h$ is also a projection. So, it suffices to show that
$$(f_1\times f_2)(g_1\times h_1,\dots,g_k\times h_k)=f_1(g_1,\dots,g_k)\times f_2(h_1,\dots,h_k)$$
for all $f_1,g_1,\dots,g_k\in{\mathcal O}_m$ and $f_2,h_1,\dots,h_k\in{\mathcal O}_n$ of suitable arities.
We have
$$(f_1\times f_2)(g_1\times h_1,\dots,g_k\times h_k)_m=(f_1\times f_2)_m((g_1\times h_1)_m,\dots,
(g_k\times h_k)_m)=$$
$$=f_1(g_1,\dots,g_k).$$
Similarly,
$$(f_1\times f_2)(g_1\times h_1,\dots,g_k\times h_k)_n=f_2(h_1,\dots,h_k).$$
Putting $f=(f_1\times f_2)(g_1\times h_1,\dots,g_k\times h_k)$,
the desired equality follows from Lemma \ref{iso0}.
\endproof

Let $\Sub({\mathcal O})$,  $\Sub({\mathcal O}_m),$ and  $\Sub({\mathcal O}_n)$ denote the lattices of
all subclones of ${\mathcal O}$, ${\mathcal O}_m$ and ${\mathcal O}_n$, respectively.

\begin{theorem} \label{iso1} 
Let  $C\in \Sub({\mathcal O})$, $D\in\Sub({\mathcal O}_m)$,  $E\in\Sub({\mathcal O}_n)$. Then
$(D\times E)_m=D$, $(D\times E)_n=E$ , and $C\subseteq C_m\times C_n$. If $C$ contains the (binary)
addition modulo $mn$ then $C_m\times C_n=C$.
\end{theorem}

\proof We have $(D\times E)_m=\{(g\times h)_m\mid g\in D,\ h\in E\ \text{of the same arity}\}=
\{g\mid g\in D\}=D$. Similarly, $(D\times E)_n=E$. If $f\in C$, then $f=f_m\times f_n\in C_m\times C_n$.

Suppose now that $C$ contains the addition modulo $mn$. Every element of $C_m\times C_n$ has
the form $g_m\times h_n$ for some $k$-ary $g,h\in C$. Since $m$ and $n$ are relatively prime,
there are  $a\in{\mathbb Z}_m$, $b\in{\mathbb Z}_n$
such that $an\equiv 1 (\od m)$ and $bm\equiv 1 (\od n)$. The~operation $f=ang+bmh=g+\dots+g+h+\dots+h$
($an$ occurrences of $g$ and $bm$ occurrences of $h$) belongs to $C$ and we claim that $f=g_m\times h_n$.
For every $x_1,\dots, x_k\in {\mathbb Z}_{mn}$ we have
$$f(x_1,\dots,x_k)\equiv g(x_1,\dots,x_k)\quad (\od m).$$
Since $g$ preserves the congruence $\od m$, we also have
$$g(x_1,\dots,x_k) \od m= g(x_1 \od m,\dots,x_k \od m) \od m=$$
$$=g_m(x_1 \od m,\dots,x_k \od m),$$
hence 
$$ f(x_1,\dots,x_k)\equiv g_m(x_1 \od m,\dots,x_k \od m)\quad (\od m),$$
and similarly,
$$ f(x_1,\dots,x_k)\equiv h_n(x_1 \od n,\dots,x_k \od n)\quad (\od n),$$
thus $f=g_m\times h_n$.
\endproof

\begin{theorem} Let $J\subseteq \Sub({\mathcal O})$ consist of all clones containing the addition modulo $mn$.
Similarly, let $J_m\subseteq \Sub({\mathcal O}_m)$ and $J_n\subseteq \Sub({\mathcal O}_n)$ consist of all clones
containing the additions modulo $m$ and modulo $n$, respectively. Then the lattice $J$ is isomorphic to the direct product
$J_m\times J_n$.\label{iso2}
\end{theorem}
\proof Clearly, if $f$ is the addition modulo $mn$, then $f_m$ and $f_n$ are additions modulo $m$ and modulo $n$, respectively. By Theorem \ref{iso1}, the assignments $C\mapsto (C_m,C_n)$ and $(D,E)\mapsto D\times E$ are mutually inverse maps between $J$ and $J_m\times J_n$. Obviously, they are order-preserving, so they
are lattice isomorphisms.
\endproof

It is also easy to see that if $f$ is the multiplication modulo $mn$, then $f_m$ and $f_n$ are the multiplications
modulo $m$ and $n$, respectively. Further, if $f$ is a (unary) constant operation, then also $f_m$ and $f_n$
are constant. It follows that if $C$ is the clone of~all polynomial operations of the group (ring) ${\mathbb Z}_{mn}$,
then $C_m$ and $C_n$ are the clones of~all polynomial operations of the groups (rings) ${\mathbb Z}_m$ and
${\mathbb Z}_n$, respectively. Hence, the isomorphism in the proof of Theorem \ref{iso2} restricts to the isomorphism
between
$\mathcal{J}_{mn}$ and $\mathcal{J}_m \times \mathcal{J}_n,$
proving Lemma \ref{LJPk}.

\section{The clone $M_1$}

From now on we consider the case of $\Zosem$.
We work with $n$-ary (ring) polynomial operations in $\Zosem$ in the form 
given in (\ref{naryf}).

Different polynomials might represent the same operation. Especially, it is easy to see that
$$0=(x-1)x(x+1)(x+2)=x^4+2x^3-x^2-2x$$
and 
$$0=(x^2-x)(y^2-y)(z^2-z)$$
for every $x,y,z\in\Zosem$,
so we can assume in (\ref{naryf}) that every $\alpha$ has at most two components greater than $1$, and no component is greater than $3$.

A polynomial of the form (\ref{naryf}) is called \emph{fully divisible} if it is divisible by $x_1x_2\dots x_n$. 
Hence, a polynomial is fully divisible iff all its monomials contain all the variables. 
A polynomial operation is called
fully divisible if it can be expressed by a~fully divisible polynomial. The importance of fully divisible polynomials lies in the following fact.
It was first observed by Bulatov for ${\mathbb Z}_{p^2}$, but the easy argument holds for $\Zosem$ too.

\begin{lemma} Let $C$ be a subclone of $P(\Zosem,+,\cdot)$ containing the addition and constants. Then $C$ is generated by its fully divisible
members.\label{fd}
\end{lemma} 
\proof
Every $f(x_1,\dots,x_n)\in C$ is a sum of fully divisible polynomial operations (with different sets of variables). We just need to show that these summands
belong to $C$. Let $f_0=f(0,x_2,\dots,x_n)$, $f_1=f-f_0$. Then $f=f_0+f_1$, $f_0,f_1\in C$. The polynomial  $f_0$ consists of the monomials of $f$ that do not
contain $x_1$, while all monomials of $f_1$ do contain $x_1$. In the next step, we split the monomials of $f_0$ and $f_1$ according to containing $x_2$.
By repeating this process, we express $f$ as a sum of fully divisible polynomials, which belong to $C$.
\endproof

One of the summands in the previous proof consists of all fully divisible monomials of $f$
(containing all variables). So we obtain the following assertion.

\begin{lemma}\label{monom} 
Let $C$ be a subclone of $P(\Zosem,+,\cdot)$ containing the addition and constants.
  Let $f\in C$, let $g$ be the sum of fully divisible monomials of $f$. 
 Then $g \in C$.
\end{lemma}

Now we are going to define the clone $M_1\subseteq P(\Zosem,+,\cdot)$.
We use the denotation $\bu=(u_1,\dots,u_k)$ for $k$-tuples of elements. Any $n$-ary operation $f:\ A^n\to A$ can be
pointwise extended to an operation $(A^k)^n\to A^k$, which will also be denoted by $f$. We say that an $n$-ary operation
$f:\ A^n\to A$ preserves a relation $R\subseteq A^k$, if $\bu^{(1)},\dots,\bu^{(n)}\in R$ implies
$f(\bu^{(1)},\dots,\bu^{(n)})\in R$.

We find it convenient to work with  $2^k$-ary relations on $\Zosem$. Our $2^k$-tuples will be indexed by $P_k$, the powerset of $\{1,\dots,k\}$.
So, we write $\bx=(x_A\mid A\in P_k)$. The components $x_A$ can be regarded as the coordinates of $\bx$ in the canonical basis. However,
in our theorems, another basis is more convenient.
For every $A\in P_k$ define $\bg^A=(g^A_B\mid B\in P_k)\in\Zosem^{P_k}$ by
$$g^A_B=\left\{ \begin{array}{l@{\quad}c}
1, \text{ if } A\subseteq B\\
0, \text{ otherwise. } \end{array} \right.$$
Notice that $\bg^\emptyset=\bjedna$, the constant $1$.
In the sequel we write $\bg^m$ instead of $\bg^{\{m\}}$ ($m=1,\dots, k$) and $\bg^0$ instead of $\bg^\emptyset$.

\begin{lemma} \label{L1.1}
Every $\bx\in\Zosem^{P_k}$ can be expressed in the form
$$\bx=\sum_{A\in P_k}a_A\bg^A,$$
where
$$a_A=(-1)^{|A|}\sum_{B\subseteq A}(-1)^{|B|}x_B.$$
for every $A\in P_k$. The expression of $\bx$ in this form is unique.
\label{ga}
\end{lemma}

\proof Let  $\bu=\sum_{A\in P_k}a_A\bg^A$. 
Clearly, $u_A=\sum_{B\subseteq A}a_B$ for every $A\in P_k$. Then $$u_A=\sum_{B\subseteq A}(-1)^{|B|}(\sum_{C\subseteq B}(-1)^{|C|}x_C)
=\sum_{C\subseteq A}(-1)^{|C|}(\sum_{C\subseteq B\subseteq A}(-1)^{|B|})x_C=x_A,$$
because for $C\subsetneq A$ we have $\sum_{C\subseteq B\subseteq A}(-1)^{|B|}=0$.

To show the uniqueness, suppose that we have two distinct expressions
$$\bx=\sum_{A\in P_k}a_A\bg^A=\sum_{A\in P_k}b_A\bg^A.$$
Then there is  a minimal set $A$ with $a_A\ne b_A$ and then
$0=x_A-x_A=\sum_{B\subseteq A}a_B-\sum_{B\subseteq A}b_B=a_A-b_A\ne 0$ a contradiction.
\endproof

The transformation of set functions using the functions $\bg^A$ is widely used in combinatorics
and elsewhere, and it is called the Möbius transform. (See \cite{grabisch2016} or \cite{hutnik}.)

The following statement is easy and will be frequently used without a reference.

\begin{lemma}For every $A,C\in P_k$,  $\bg^A\cdot \bg^C=\bg^{A\cup C}$ under the pointwise
multiplication.
\end{lemma}

Let $Z$ be the $2^4$-ary relation on $\Zosem$ consisting of all tuples $\bu=(u_A\mid A\in P_4)$ satisfying
\begin{enumerate}
\item[(Z1)] $a_{\{2\}}\equiv 2a_{\{1\}}\ (\od 4)$, 
$a_{\{4\}}\equiv 2a_{\{3\}}\ (\od 4)$;
\item[(Z2)] $a_A\equiv 0\ (\od 2)$ whenever $|A|\ge 2$;
\item[(Z3)] $a_A\equiv 0\ (\od 4)$ whenever $|A|\ge 2$, $A\cap\{2,4\}\ne\emptyset$;
\item[(Z4)] $a_A=0$ whenever $\{2,4\}\subseteq A$.
\end{enumerate}

Let $M_1$ be the clone of all ring polynomials preserving
the relation $Z$.

Every clone of ring polynomials containing the addition and constants is generated by its fully divisible members. So, our first aim is to characterize the fully divisible members of $M_1$. 
As shown in the previous section, the general form of a~fully divisible $n$-ary polynomial operation on $\Zosem$ is
\begin{equation}f(\bx)=x_1x_2\dots x_n\sum_{\alpha\in I_n}b_\alpha\bx^\alpha,\label{gen}\end{equation}
where $I_n$ is the set of all $n$-tuples from $\{0,1,2\}^n$ with at most two nonzero components.

It is easy to see that the addition preserves $Z$.

\begin{lemma} The operation $f(x)=2x_1x_2\dots x_n$  preserves $Z$ for every $n\ge 1$.
\end{lemma}

\proof
We need to show that  $f(\bu^{(1)},\dots,\bu^{(n)})\in Z$ whenever $\bu^{(i)}\in Z$ for every $i$. 
Consider $\bu^{(i)}$ in the form 
$$\bu^{(i)}=\sum_{A\in {\mathcal P}(4)}a_A^{(i)}\bg^A,$$
where $a_A^{(i)}=0$ for $\{2,4\}\subseteq A$. Then 
$$f(\bu^{(1)},\dots,\bu^{(n)})=\sum_{A\in {\mathcal P}(4)}a_A\bg^A,$$
where 
$$a_A=\sum \{2a_{A_1}^{(1)}\cdot a_{A_2}^{(2)}\cdot\dots a_{A_n}^{(n)}\mid A_1\cup\dots\cup A_n=A\}.$$

Clearly, every $a_A$ is even, so (Z2) is satisfied. By (Z1), all $a_{\{2\}}^{(i)}$ and $a_{\{4\}}^{(i)}$ are even. Consequently, $a_{\{2\}}$ and $a_{\{4\}}$ are
divisible by $4$, so $f(\bu^{(1)},\dots,\bu^{(n)})$ satisfies (Z1). Since every $a_A^{(i)}$ with $A$ containing $2$ or $4$ is even, we obtain (Z3). Finally,
if $\{2,4\}\subseteq A$, then  every summand in the expression for $a_A$  contains $2a_{\{2\}}^{(i)}a_{\{4\}}^{(j)}$ for some $i,j$, hence $a_A=0$,
so (Z4) holds.
\endproof

Every polynomial operation with even coefficients belongs to the clone generated by the addition and all operations $2x_1\dots x_n$. Hence, we obtain the following.

\begin{lemma} Every operation of the form (\ref{gen}) with all $b_\alpha$ even preserves $Z$.\label{suffices}
\end{lemma}

Now we prove the converse. The unary case needs to be treated separately.

\begin{lemma} The operation $f(x)=ax^3+bx^2+cx$ preserves $Z$ if and only if both $a$ and $b$ are even.\label{unary}\end{lemma}
\proof All linear operations obviously preserve $Z$. Hence, $f$ preserves $Z$ if and only if $g(x)=ax^3+bx^2$
preserves it.

First, suppose that $g$ preserves $Z$. Let $\bu=\bg^{\{1\}}+2\bg^{\{2\}}$. Obviously, $\bu\in Z$, so $g(\bu)\in Z$.
It is easy to calculate that $g(\bu)=(a+b)\bg^{\{1\}}+4b\bg^{\{2\}}+(2a+4b)\bg^{\{1,2\}}$. From (Z3) we obtain
$2a+4b\equiv 0 (\mod 4)$, so $a$ is even. From (Z1) we obtain $2(a+b)\equiv 4b\equiv 0\ (\mod 4)$, so $b$ is even.

The converse follows from Lemma \ref{suffices}.
\endproof

\begin{lemma} Let $f$ of the form (\ref{gen}) preserve $Z$. Suppose that $n\ge 2$, $b_\alpha\in\{0,1\}$ for every $\alpha$. 
Then $b_\alpha=0$ for every $\alpha$.\label{01}
\end{lemma}
\proof  Let $\bu^{(1)}=\bg^{\{1\}}+2\bg^{\{2\}}$, $\bu^{(2)}=\bg^{\{3\}}+2\bg^{\{4\}}$ and 
$\bu^{(i)}=\bg^\emptyset=\bjedna$ (the constant tuple) for $i=3,\dots,n$. Clearly, $\bu^{(i)}\in Z$ for every $i$,
hence $f(\bu^{(1)},\dots,\bu^{(n)})\in Z$. Let $g(x_1,x_2)=f(x_1,x_2,1,1,\dots,1)$. Then $g$ is a binary polynomial, which can be expressed 
in the form
$$g(\bx)=x_1x_2\sum_{\alpha\in I_2}c_\alpha\bx^\alpha.$$
We have $g(\bu^{(1)},\bu^{(2)})=f(\bu^{(1)},\dots,\bu^{(n)})\in Z$.
Clearly, 
$$c_{(\alpha_1,\alpha_2)}=b_{(\alpha_1,\alpha_2,0,\dots,0)}$$
whenever $\alpha_1,\alpha_2\ne 0$. 

Expressing $g(\bu^{(1)},\bu^{(2)})$ in the form $\sum_{A\in P_n}d_A\bg^A$ requires a fair amount of calculations.
We have $(\bu^{(1)})^2=\bg^{\{1\}}+4\bg^{\{1,2\}}+4\bg^{\{2\}}$,  $(\bu^{(1)})^3=\bg^{\{1\}}+2\bg^{\{1,2\}}$, 
$(\bu^{(2)})^2=\bg^{\{3\}}+4\bg^{\{3,4\}}+4\bg^{\{4\}}$,  $(\bu^{(2)})^3=\bg^{\{3\}}+2\bg^{\{3,4\}}$.
Hence,
\begin{equation}\label{e33}
(\bu^{(1)})^3(\bu^{(2)})^3=(\bg^{\{1\}}+2\bg^{\{1,2\}})(\bg^{\{3\}}+2\bg^{\{3,4\}})=\end{equation}
$$= \bg^{\{1,3\}}+2\bg^{\{1,3,4\}}+2\bg^{\{1,2,3\}}+4\bg^{\{1,2,3,4\}}.$$
Similarly we find
\begin{equation}
(\bu^{(1)})^3(\bu^{(2)})^2= \bg^{\{1,3\}}+4\bg^{\{1,3,4\}}+4\bg^{\{1,4\}}+2\bg^{\{1,2,3\}},\end{equation}
\begin{equation}
(\bu^{(1)})^2(\bu^{(2)})^3= \bg^{\{1,3\}}+4\bg^{\{1,2,3\}}+4\bg^{\{2,3\}}+2\bg^{\{1,3,4\}},\end{equation}
\begin{equation}\label{e22}
(\bu^{(1)})^2(\bu^{(2)})^2= \bg^{\{1,3\}}+4\bg^{\{1,3,4\}}+4\bg^{\{1,2,3\}}+4\bg^{\{1,4\}}+4\bg^{\{2,3\}},\end{equation}
\begin{equation}\label{e31}
(\bu^{(1)})^3\bu^{(2)}= \bg^{\{1,3\}}+4\bg^{\{1,2,4\}}+2\bg^{\{1,2,3\}}+2\bg^{\{1,4\}},\end{equation}
\begin{equation}
\bu^{(1)}(\bu^{(2)})^3= \bg^{\{1,3\}}+4\bg^{\{2,3,4\}}+2\bg^{\{1,3,4\}}+2\bg^{\{2,3\}},\end{equation}
\begin{equation}\label{e21}
(\bu^{(1)})^2\bu^{(2)}= \bg^{\{1,3\}}+4\bg^{\{1,2,3\}}+4\bg^{\{2,3\}}+2\bg^{\{1,4\}},\end{equation}
\begin{equation}\label{e12}
\bu^{(1)}(\bu^{(2)})^2= \bg^{\{1,3\}}+4\bg^{\{1,3,4\}}+4\bg^{\{1,4\}}+2\bg^{\{2,3\}},\end{equation}
\begin{equation}\label{e11}
\bu^{(1)}\bu^{(2)}= \bg^{\{1,3\}}+4\bg^{\{2,4\}}+2\bg^{\{1,4\}}+2\bg^{\{2,3\}}.\end{equation}

Now we have
$$g(\bu^{(1)},\bu^{(2)})=\sum_{A\in P_n}d_A\bg^A=\sum_{i,j=1}^3c_{(i-1,j-1)}(\bu^{(1)})^i(\bu^{(2)})^j.$$
The tuple $\bg^{\{1,2,3,4\}}$ only appears in (\ref{e33}), so $d_{\{1,2,3,4\}}=4c_{(2,2)}$. By the condition (Z4)
we have $4c_{(2,2)}=d_{\{1,2,3,4\}}=0$, so  $c_{(2,2)}$ is even. Since $c_{(2,2)}=b_{(2,2,0,\dots,0)}\in\{0,1\}$,
we obtain $b_{(2,2,0,\dots,0)}=0$. By symmetry, $b_\alpha=0$ whenever $\alpha$ has any two components equal to $2$.
The tuple $\bg^{\{1,2,4\}}$ only appears in (\ref{e31}), so $d_{\{1,2,4\}}=4c_{(2,0)}$. By the condition (Z4)
we have $4c_{(2,0)}=d_{\{1,2,4\}}=0$, so  $c_{(2,0)}$ is even. (We cannot yet claim $c_{(2,0)}=0$.)
The tuple $\bg^{\{1,2,3\}}$ appears in multiple places, so we obtain
$$d_{\{1,2,3\}}=2c_{(2,2)}+2c_{(2,1)}+4c_{(1,2)}+4c_{(1,1)}+2c_{(2,0)}+4c_{(1,0)}.$$
By the condition (Z3), $d_{\{1,2,3\}}\equiv 0\ (\mod 4)$. Since $c_{(2,2)}=0$ and $c_{(2,0)}$ is even,
we obtain that $c_{(2,1)}$ is even. Since $c_{(2,1)}=b_{(2,1,0,\dots,0)}\in\{0,1\}$, we have
 $c_{(2,1)}=b_{(2,1,0,\dots,0)}=0$. By the symmetry, $b_\alpha=0$ whenever $\alpha$ contains both $2$ and $1$.

Next,
$$c_{(2,0)}=\sum\{b_\alpha\mid \alpha_1=2,\ \alpha_2=0\}.$$
We have already proved that all the summands are zero, except $b_{(2,0,0,\dots,0)}$. Hence,
$c_{(2,0)}=b_{(2,0,\dots,0)}\in\{0,1\}$. Since $c_{(2,0)}$ is even, we obtain $b_{(2,0,\dots,0)}=0$. By symmetry,
$b_\alpha=0$ whenever $\alpha$ contains $2$. 

Up to now, we have $c_{(2,2)}=c_{(2,1)}=c_{(1,2)}=c_{(2,0)}=c_{(0,2)}=0$. Further, $\bg^{\{2,4\}}$ appears only in (\ref{e11}), so $d_{\{2,4\}}=4c_{(0,0)}$. By (Z4), $c_{(0,0)}$ is even.
For $d_{\{2,3\}}$ we obtain
$$d_{\{2,3\}}=4c_{(1,1)}+4c_{(1,0)}+2c_{(0,1)}+2c_{(0,0)}.$$
Since (by (Z3)), $d_{\{2,3\}}$ is a multiple of $4$, and $c_{(0,0)}$ is even, $c_{(0,1)}$ must be even too.
By a similar argument (using $d_{\{1,4\}}$) we show that $c_{(1,0)}$ is even. Now we have
\begin{equation}\label{n13}
d_{\{1,3\}}=c_{(1,1)}+c_{(1,0)}+c_{(0,1)}+c_{(0,0)}.\end{equation}
Since (by (Z2)) $d_{\{1,3\}}$ is even, $c_{(1,1)}$ is even too. However, $c_{(1,1)}=b_{(1,1,0,\dots,0)}\in\{0,1\}$,
hence $b_{(1,1,0,\dots,0)}=0$. By symmetry, $b_\alpha=0$ whenever $\alpha$ has two components equal to $1$.

Now we consider
$$c_{(1,0)}=\sum\{b_\alpha\mid \alpha_1=1,\ \alpha_2=0\}.$$
We have already proved that all the summands are zero, except $b_{(1,0,0,\dots,0)}$. Hence,
$c_{(1,0)}=b_{(1,0,\dots,0)}\in\{0,1\}$. Since $c_{(1,0)}$ is even, we obtain $b_{(1,0,\dots,0)}=0$. By symmetry,
$b_\alpha=0$ whenever $\alpha$ contains $1$. 

It remains to show that $b_{(0,0,\dots,0)}=0$. We have
$$c_{(0,0)}=\sum\{b_\alpha\mid \alpha_1=0,\ \alpha_2=0\}=b_{(0,0,\dots,0)},$$
as all other summands are zero. By our assumption, $b_{(0,0,\dots,0)}\in\{0,1\}$. Hence,
$b_{(0,0,\dots,0)}=0$, which completes the proof.
\endproof

\begin{theorem}Let $f$ be defined by (\ref{gen}) with $n\ge 2$. Then $f$ preserves $Z$ if and only
 if every $b_\alpha$ is even.\label{M}
\end{theorem}

\proof Every $f$ of the form (\ref{gen}) can be expressed as the sum
\begin{equation}\label{decompose}
f(\bx)=x_1x_2\dots x_n\sum_{\alpha\in I_n}c_\alpha\bx^\alpha + x_1x_2\dots x_n\sum_{\alpha\in I_n}k_\alpha\bx^\alpha,\end{equation}
where every $c_\alpha$ is even and every $k_\alpha$ belongs to $\{0,1\}$.
Suppose that $f$ preserves $Z$. By Lemma \ref{suffices}, the first summand in (\ref{decompose}) also preserves $Z$. 
The addition also preserves  $Z$, so $ x_1x_2\dots x_n\sum_{\alpha\in I_n}k_\alpha\bx^\alpha$ preserves $Z$.
By Lemma \ref{01}, the second summand in (\ref{decompose}) is zero, so $f(\bx)=x_1x_2\dots x_n\sum_{\alpha\in I_n}c_\alpha\bx^\alpha$.

The inverse implication follows from Lemma \ref{suffices}.
\endproof

So, we study fully divisible polynomials of the form
\begin{equation}f=2x_1\dots x_n(\sum_{\alpha\in I_n}b_\alpha\bx^\alpha),\label{fold}.\end{equation}
 This form can be further simplified. Because of the equality
$$2(x^2+x)(y^2+y)=0$$
for every $x,y\in\Zosem$, we can replace all monomials $b_\alpha\bx^\alpha$, where $\alpha$ has two nonzero components,
with monomials of lower degrees. We obtain that the operation defined by (\ref{fold}) is equal to the operation defined
by a polynomial of the form
\begin{equation}f=2x_1\dots x_n(\sum_{i=1}^na_ix_i^2+\sum_{i=1}^nb_ix_i+c).\label{f}\end{equation}

The equality $4x^3=4x^2=4x$, which holds for all $x\in\Zosem$ means that $f$ of the~form (\ref{f}) does not change if we simultaneously replace some $a_k$ or $b_k$ with
$a_k-2$ (or $b_k-2$, respectively) and $c$ with $c+2$. Also, $f$ does not change if we replace $c$ with $c+4$.  Hence, from now on, we will assume that in (\ref{f})
$$a_i,b_i\in\{0,1\}\ \text{for every} \ i\ \text{and}\ c\in\{0,1,2,3\}.$$

In the sequel, let ${\mathcal G}_1$ be the set of all polynomials of this form (for every $n\ge 1$). Further, let ${\mathcal G}_0$ be the set of all polynomials of the form (\ref{f})
such that $\sum_{i=1}^na_i+\sum_{i=1}^nb_i+c$ is even.

\begin{lemma} Let $f$ be an operation of the form (\ref{f}) with $n\ge 2$. Then $f\in{\mathcal G}_0$ if and only if $f$
preserves the relation
$$M=\{\bx\in Z\mid a_{\{1,3\}}\equiv 0\ (\od 4)\}.$$
\end{lemma}
\proof Suppose that $f\in{\mathcal G}_0$, $\bu^{(1)},\dots,\bu^{(n)}\in M$. Since $M\subseteq Z$ and $f$ preserves $Z$, we obtain
$f(\bu^{(1)},\dots,\bu^{(n)})\in Z$. The additional condition for $M$ is satisfied, because clearly $\rng(f)\subseteq\{0,4\}$.

Conversely, suppose that $f\in{\mathcal G}_1\setminus{\mathcal G}_0$. Then $g=f+2x_1\dots x_n\in{\mathcal G}_0$. By the first part of this proof,
$g$ preserves $M$. Hence, $f$ preserves $M$ if and only if $2x_1\dots x_n$ does. Let
$\bu^{(1)}=\bg^{\{1\}}+2\bg^{\{2\}}$, $\bu^{(2)}=\bg^{\{3\}}+2\bg^{\{4\}}$ and 
$\bu^{(i)}=\bg^\emptyset=\bjedna$ (the constant tuple) for $i=3,\dots,n$. Clearly, $\bu^{(i)}\in M$ for every $i$, and
$$2\bu^{(1)}\dots\bu^{(n)}=2\bg^{\{1,3\}}+4\bg^{\{1,4\}}+4\bg^{\{2,3\}}\notin M,$$
so $2x_1\dots x_n$ and $f$ do not preserve $M$.
\endproof

\section{The subclones of $M_1$}

For every $f\in{\mathcal G}_1$, let $C(f)$ denote the clone on $\Zosem$
generated by $f$ together with the addition and the (unary) constant operations.
Now we investigate equalities and inequalities between these clones.

We introduce the following denotations for some $n$-ary ($n\ge 1$) polynomials over ${\mathbb Z}_8$.
$$r_n=x_1x_2\dots x_n;$$
$$t_n=\left\{ \begin{array}{l@{\quad}c}
x_1x_2\dots x_n(x_1+\dots +x_n), \text{ if } n\text{ is even}\\
x_1x_2\dots x_n(x_1+\dots +x_n+1), \text{ if } n\text{ is odd.}
\end{array} \right.$$
$$s_n=\left\{ \begin{array}{l@{\quad}c}
x_1x_2\dots x_n(x_1+\dots +x_n), \text{ if } n\text{ is odd}\\
x_1x_2\dots x_n(x_1+\dots +x_n+1), \text{ if } n\text{ is even.}
\end{array} \right.$$
$$p_n=x_1^2x_2\dots x_n;$$
$$q_n=x_1^3x_2\dots x_n;$$
$$w_n=x_1x_2\dots x_n(x_1^2+1);$$
$$u_n=x_1x_2\dots x_n(x_1+1);$$
$$v_n=x_1x_2\dots x_n(x_1+x_2).$$

(The operation $v_1$ is undefined.) 
The operations $2r_n$, $4r_n$, $2t_n$, $2s_n$, $2p_n$, $2q_n$, $2w_n$, $2u_n$ and $2v_n$ belong to ${\mathcal G}_1$. Out of them, $4r_n$, $2t_n$, $2w_n$, $2u_n$ and $2v_n$ belong to ${\mathcal G}_0$. Notice that  $4t_n$, $4w_n$, $4u_n$ and $4v_n$
are zero.

\begin{lemma}\label{L3.1} Let $f\in{\mathcal G}_1$ be of the form (\ref{f}).
\begin{enumerate}
\item[(i)] If $a_i=1$ for some $i$, then $4r_k\in C(f)$ for every $k\le n+2$;
\item[(ii)] If $b_i=1$ for some $i$, then $4r_k\in C(f)$ for every $k\le n+1$;
\end{enumerate}
\end{lemma}

\proof (i) We can assume that $a_1=1$ and write $f$ in the form
$$f=2x_1^3x_2\dots x_n+2x_1\dots x_n(\sum_{i=2}^na_ix_i^2+\sum_{i=1}^nb_ix_i+c).$$
Now we make a substitution $x_1:=x_1+x_{n+1}+x_{n+2}$ (other variables unchanged).
After this substitution, we obtain a polynomial operation $g$ that belongs to $C(f)$. The polynomial $g$
contains a single monomial divisible by  $x_1x_2\dots x_{n+2}$ and it is
$$2\cdot (6x_1x_{n+1}x_{n+2})\cdot x_2\dots x_n=4r_{n+2}.$$
By Lemma \ref{monom} we obtain $4r_{n+2}\in C(f)$. Substituting $x_i=1$ for all $i>k$ we obtain $4r_k\in C(f)$
for every $k\le n+2$.

(ii) Because of (i), we can assume that $a_i=0$ for all $i$ and $b_1=1$. Then
$$f=2x_1^2x_2\dots x_n+2x_1\dots x_n(\sum_{i=2}^nb_ix_i+c).$$
After the substitution $x_1:=x_1+x_{n+1}$ we obtain a polynomial operation $g\in C(f)$, which has a single monomial
divisible by $x_1\dots x_{n+1}$, namely
$$2\cdot(2x_1x_{n+1})x_2\dots x_n=4r_{n+1}.$$
Hence, $4r_k\in C(f)$ for every $k\le n+1$.

\endproof

Hence, if $f$ is any non-zero member of ${\mathcal G}_1$, then $C(f)$ contains $4r_n$. Consequently,
if $f$ has at least one $a_i$ or $b_i$ nonzero, then $C(f)=G(g)$ for
$$g=2x_1\dots x_n(\sum_{i=1}^na_ix_i^2+\sum_{i=1}^nb_ix_i+c-2).$$ 
Hence, we only need $c=2$ to include the operations $4r_n$ in our considerations. In all other cases, we can assume $c\in\{0,1\}$.

\begin{lemma} For every $n\ge 1$,
$2r_n\in C(2r_{n+1})$, $4r_n\in C(4r_{n+1})$,$2t_n\in C(2t_{n+1})$, $2s_n\in C(2s_{n+1})$,$2p_n\in C(2p_{n+1})$,
$2q_n\in C(2q_{n+1})$, $2u_n\in C(2u_{n+1})$,
and 
$2w_n\in C(2w_{n+1})$. For every $n\ge 2$, $2v_n\in C(2v_{n+1})$.\label{chains}
\end{lemma}
\proof All the statements can be shown by substituting
$x_{n+1}=1$. Only the cases of $2t_n$ for even $n$ and
$2s_n$ for odd $n$ require one more argument. In these cases
we have $2t_{n+1}(x_1,\dots,x_n,1)=2t_n+4r_n$ and
$2s_{n+1}(x_1,\dots,x_n,1)=2s_n+4r_n$, so 
$2t_n\in C(2t_{n+1})$ and $2s_n\in C(2s_{n+1})$ follow
from Lemma \ref{L3.1}.
\endproof

\begin{lemma} Let  $f\in{\mathcal G}_1$ be given by (\ref{f}). Assume that $a_i\ne 0$ for some $i$. Then $2u_n\in C(f)$.\label{un}
\end{lemma}

\proof We can assume $a_1=1$ and write $f=2x_1\dots x_n(x_1^2+b_1x_1+h(x_2,\dots,x_n))$.
We calculate
$$f(x_1+1,x_2,\dots,x_n)-f(x_1,\dots,x_n)=$$
$$=2x_1\dots x_n((x_1+1)^2+b_1(x_1+1)+h(x_2,\dots,x_n))+$$
$$+2x_2\dots x_n((x_1+1)^2+b_1(x_1+1)+h(x_2,\dots,x_n))+$$
$$-2x_1\dots x_n(x_1^2+b_1x_1+h(x_2,\dots,x_n))=$$
$$=2x_1\dots x_n(2x_1+1+b_1)+2x_2\dots x_n((x_1+1)^2+b_1(x_1+1)+h(x_2,\dots,x_n)).$$

Omitting the polynomial $2x_2\dots x_n(1+b_1+h(x_2,\dots,x_n))$ which does not contain $x_1$, we see that
$C(f)$ contains
$$2x_1\dots x_n(2x_1+1+b_1+x_1+2+b_1)=2u_n+2x_1\dots x_n(2x_1+2b_1+2).$$
The summand $2x_1\dots x_n(2x_1+2b_1+2)=4x_1\dots x_n(x_1+b_1+1)$ is zero or $4r_n$, so it belongs to
$C(f)$. Then $2u_n\in C(f)$, too.
\endproof

\begin{lemma}\label{un3}
Let $f\in{\mathcal G}_1$ be given by (\ref{f}). 
Assume that $a_i\ne 0$ for some $i$. 
Then $2w_n\in C(f)$ or $2q_n\in C(f)$.
\end{lemma}

\proof 
We can assume $a_1=1$. First we show that $C(f)$ contains $g\in{\mathcal G}_1$ of the form
\begin{equation}g=2x_1\dots x_n(x_1^2+\sum_{i=1}^nb_i'x_i+c')\label{g}.\end{equation}
We need to eliminate all nonzero $a_i$'s except $a_1$. Let us suppose that $a_2=1$. 
Then we substitute $x_1:=x_1+x_2$ in $f$ with other variables unchanged. 
We obtain a~member of $C(f)$, which is the sum of
$$f_1=2x_1\dots x_n((x_1+x_2)^2+x_2^2+\sum_{i=3}^na_ix_i^2+b_1x_2+\sum_{i=1}^nb_ix_i+c)$$
and 
$$f_2=2x_2^2x_3\dots x_n((x_1+x_2)^2+x_2^2+\sum_{i=3}^na_ix_i^2+b_1x_2+\sum_{i=1}^nb_ix_i+c).$$
Some monomials of $f_2$ do not contain $x_1$. After omitting them, we obtain the~polynomial
$$f_3=2x_2^2x_3\dots x_n(x_1^2+2x_1x_2+b_1x_1)=2x_1\dots x_n(x_1x_2+2x_2^2+b_1x_2),$$
and we still have $f_1+f_3\in C(f)$. Using the identity
$$2x_1^2x_2^2=-2x_1^2x_2-2x_1x_2^2-2x_1x_2$$
we have
$$f_3=2x_1\dots x_n(-x_1-x_2-1+2x_2^2+b_1x_2).$$
Further, the monomials $2x_1\dots x_n\cdot2x_1x_2$ and $2x_1\dots x_n\cdot 2x_2^2$ are both equal to $4r_n$. 
So we have
$$f_1=2x_1\dots x_n(x_1^2+2+2+\sum_{i=3}^na_ix_i^2+b_1x_2+\sum_{i=1}^nb_ix_i+c),$$
$$f_3=2x_1\dots x_n(-x_1-x_2-1+2+b_1x_2).$$
The sum $g_1=f_1+f_3$ belongs to $C(f)$ and can be written in the form
$$g_1=2x_1\dots x_n(x_1^2+\sum_{i=3}^na_ix_i^2+\sum_{i=1}^nb_i'x_i+c').$$
Repeating this process, we obtain $g\in C(f)$ of the required form (\ref{g}).

By Lemma \ref{un}, $2u_{n,i}=2x_1\dots x_n(x_i+1)$ belongs to $C(f)$ for every $i$. Hence,
$$h=g-\sum_{i=1}^n2b_i'u_{n,i}=2x_1\dots x_n(x_1^2+c'-\sum_{i=1}^nb_i')=2w_n+2dr_n=2q_n+2(d+1)r_n,$$
where $d=c'-\sum_{i=1}^nb_i'-1$, also belongs to $C(f)$. If $d$ is even, then $2dr_n\in C(f)$ and hence also $2w_n\in C(f)$. If $d$ is odd, then $2(d+1)r_n\in C(f)$ and $2q_n\in C(f)$.
 \endproof

\begin{lemma}\label{L3.6}
   For every $n \geq 2$, $C(2w_n)=C(2v_{n+1})$ holds .
\end{lemma}
\proof
    Firstly, we will show that we can obtain $2w_n$ from the operation $2v_{n+1}$. Since $n\ge 2$, we can substitute
    $x_{n+1}:=x_1$.
    We have
           $$2v_{n+1}(x_1, x_2, \dots, x_{n}, x_1)=2x_1^3 x_2 \dots x_{n} + 2x_1^2 x_2^2 \dots x_{n}.$$
    The first summand is equal to $2w_n(x_1,\dots,x_n)-
    2r_n(x_1,\dots,x_n)$. For the second summand we use the equation 
    $2x_1^2x_2^2 = 6x_1^2 x_2+6x_1x_2^2 +6x_1x_2$
    and obtain
    $$2x_1^3 x_2 \dots x_{n} +6x_1^2 x_2 \dots x_{n}+6x_1 x_2^2 \dots x_{n}+6x_1 x_2 \dots x_{n}=$$
$$=2w_n (x_1, x_2, \dots, x_{n})+ 6v_n(x_1, x_2, \dots, x_{n})+4r_n (x_1, x_2, \dots, x_{n}).$$
According to Lemma \ref{L3.1} the polynomial $4r_n$ belongs to $C(2v_{n+1})$. 
 Also, $2v_{n}\in C(2v_{n+1})$, so the clone $ C(2v_{n+1})$ includes the function $ 2w_{n}$, too.

 Conversely, we prove that $2v_{n+1}\in C(2w_{n})$.
 We have
 $$2w_{n}(x_1 +x_{n+1}, x_2, \dots, x_{n})= 2(x_1 +x_{n+1}) x_2 \dots x_{n} ((x_1 +x_{n+1})^2 + 1).$$
The fully divisible part of this polynomial belongs to $C(2w_n)$ and is equal to
 $$2x_1x_2\dots x_{n+1}(x_{n+1}+2x_{n+1}+ x_1+2x_1)=2v_{n+1}+4x_1x_2\dots x_{n+1}^2+4x_1^2x_2\dots x_{n+1}.$$
Since $4x_1^2=4x_1$ and $4x_{n+1}^2=4x_{n+1}$, we have $4x_1x_2\dots x_{n+1}^2+4x_1^2x_2\dots x_{n+1} =0$,
so $2v_n \in C(2w_{n}) $.
\endproof

\begin{lemma}\label{Lemma 3.5} For every $n\ge 1$, $2w_n\in C(2q_n)$ and $2r_n\in C(2q_n)$.\label{wpodq}\end{lemma}
\proof Since $2r_1=2x_1$ is linear, we have $2r_1\in
C(2q_1)\cap C(2w_1)$ and $C(2q_1)=C(2w_1)$. Suppose now $n\ge 2$. The operation $2(x_1+x_2)^3x_3\dots x_{n+1}$ belongs to
$C(2q_n)$.
The~fully divisible part of this polynomial is $6x_1\dots x_{n+1}(x_1+x_2)=-2v_{n+1}$.
Hence, $2v_{n+1}\in C(2q_n)$. By Lemma \ref{L3.6}, $2w_n\in C(2q_n)$. Since $2r_n=2w_n-2q_n$,
we also obtain $2r_n\in C(2q_n)$.
\endproof

As a consequence, we obtain the next lemma.

\begin{lemma}
Let $f\in{\mathcal G}_1$ be given by (\ref{f}). 
Assume that $a_i\ne 0$ for some $i$. 
Then $2w_n\in C(f)$.
\end{lemma}

\begin{theorem} Let $f\in{\mathcal G}_1$ be given by (\ref{f}). Assume that $a_i\ne 0$ for some $i$. If $f\in{\mathcal G}_0$, then $C(f)=C(2w_n)$. If $f\in{\mathcal G}_1\setminus{\mathcal G}_0$,
then $C(f)=C(2q_n)$.
\end{theorem}

\proof  Let $f$ be $n$-ary,  given by (\ref{f}). We already know that, for every $i$, the~operations $2w_{n,i}=2x_1\dots x_n(x_i^2+1)$ and
$2u_{n,i}=2x_1\dots x_n(x_i+1)$ belong to $C(2w_n)$. Then also
$$g=\sum_{i=1}^n2a_iw_{n,i}+\sum_{i=1}^n2b_iu_{n,i}$$
belongs to $C(2w_n)$. Clearly,
$$g=2x_1\dots x_n(\sum_{i=1}^na_ix_i^2+\sum_{i=1}^nb_ix_i+d),$$
where $d=\sum_{i=1}^na_i+\sum_{i=1}^nb_i$.
Hence, $f=g+2x_1\dots x_n(c-d)$. Then $g\in C(2w_n)\subseteq C(2q_n)$ 
and $2r_n\in C(2q_n)$ imply $f\in C(2q_n)$. So, we have $C(2w_n)\subseteq C(f)\subseteq C(2q_n)$.

If $f\in{\mathcal G}_0$, the value $c-d$ is even, hence 
$2x_1\dots x_n(c-d)\in C(2w_n)$. Therefore, $f\in C(2w_n)$.

If $f\notin{\mathcal G}_0$, then $d-c$ is odd and $2r_n=f-g-2(d-c+1)r_n\in C(f)$.
Then also $2q_n=2w_n-2r_n\in C(f)$.
\endproof

\begin{theorem}
    Suppose that $f\in \mathcal{G}_1$, $a_i=0$ for all $i$ and $b_j=1$ for some $j$. 
    Let $K=\{j\mid \ b_j=1\}$.
    Then
    \begin{itemize}
        \item[(i)] if $K=\{1, \dots, n\}$, then $f=2t_n$ or $f=2s_n$;
        \item[(ii)] if $1\leq |K|< n$, then $2v_n\in C(f)$;
        \item[(iii)] if $1\leq |K|< n$ and $|K| $ is even, then $C(f)=C(2v_n)$ or $C(f)= C(2v_n)\vee C(2r_n)$;
        \item[(iv)] if $1\leq |K|< n$ and $|K| $ is odd, then $C(f)=C(2u_n)$ or $C(f)=C(2p_n)$.
   
    \end{itemize}
\end{theorem}

\proof
Let us denote $|K|=k$.  Without loss of generality $K=\{1, \dots, k\}$, which means that 
        $$ f= 2x_1 \dots x_n(x_1 +\dots +x_k )$$
        or
$$ f= 2x_1 \dots x_n(x_1 +\dots +x_k +1).$$ So, (i) is clear.

  (ii) We have 
    $$f(x_1, \dots, x_n)+f(x_2, \dots, x_n, x_1)=$$
    $$=2v_n(x_1, x_{k+1}, x_2, \dots, x_n) + 4x_1\dots x_n(x_2+x_3+ \dots +x_k)+4mx_1\dots x_n,$$
    where $m\in\{0,1\}$ depending on the form of $f$. 
    For every $i$ we have $4x_1\dots x_nx_i=4r_n\in C(f)$ by Lemma \ref{L3.1}.
    So, $2x_1\dots x_n(x_1+x_{k+1})\in C(f)$ and therefore $2v_n $ is included in $C(f)$, too.

(iii) If $|K|$ is even, then either
$$f=2v_n(x_1, x_2, \dots , x_n )+2v_n(x_3,x_4, \dots , x_n, x_1, x_2)+\dots+$$
$$ +2v_n(x_{k-1},x_k, \dots, x_n,\dots, x_{k-2}),$$
or
$$f=2v_n(x_1, x_2, \dots , x_n )+2v_n(x_3,x_4, \dots , x_n, x_1, x_2)+\dots+$$
$$ +2v_n(x_{k-1},x_k, \dots, x_n,\dots, x_{k-2})+2r_n.$$
In the first case we have $f\in C(2v_n)$, so (together with (ii)) $C(f)=C(2v_n)$.
In the~second case we have $f\in C(2v_n)\vee C(2r_n)$, so $C(f)\subseteq C(2v_n)\vee C(2r_n)$.
The~same equality implies $2r_n\in C(f)$ (since $2v_n\in C(f)$), so 
$C(2v_n)\vee C(2r_n)\subseteq C(f)$.

(iv) 
If $|K|$ is odd, then either
$$f=2u_n+2v_n(x_2, x_3, \dots) + 2v_n(x_4, x_5, \dots)+\dots +2v_n(x_{k-1}, x_k, \dots),$$
or
$$f=2p_n+2v_n(x_2, x_3, \dots) + 2v_n(x_4, x_5, \dots)+\dots +2v_n(x_{k-1}, x_k, \dots).$$
From (ii) we know that $2v_n\in C(f)$ and also $2v_n\in C(2u_n)$, $2v_n\in C(2p_n)$.
Hence, in the first case, we obtain $C(2u_n)= C(f)$, and in the second case $C(f)=C(p_n)$.
 \endproof

So, every clone in the interval $\langle P(\mathbb{Z}_{8}, +), M_1 \rangle$ is generated by its members $f$ of the form $2t_n$, $2u_n$, $2v_n$,  $2s_n$,
$2p_n$, $2q_n$, $2r_n$ and $4r_n$. Notice that the operations $2w_n$ are no longer needed, as $C(2w_n)=C(2v_{n+1})$
for $n\ge 2$ and $C(2w_1)=C(2q_1)$.

Since clones are generated by operations of each type form a chain with increasing $n$ (for instance, $C(2t_2)\subseteq C(2t_3)\subseteq C(2t_4)\subseteq\dots$),
we obtain the following result.

\begin{theorem}
Every clone in $\langle P(\mathbb{Z}_{8}, +), M_1 \rangle$ can be expressed in the form of a join
$$C(2t_{n_1})\vee C(2u_{n_2})\vee C(2v_{n_3})\vee C(2s_{n_4})\vee C(2p_{n_5})\vee C(2q_{n_6})\vee C(2r_{n_7})\vee C(4r_{n_8}).$$
In this join we admit $n_i=0$ (and $n_3=1$), which means that the corresponding joinand is not present, and $n_i=\infty$, which means that the clone contains all operations
of the corresponding type.
\end{theorem}
 Of course, some joinands might be redundant, so an expression of a clone in the above form is not unique. In the forthcoming
 continuation of this paper we provide an algorithm for
 comparing two clones of this form, which completes
 the description of the interval  $\langle P(\mathbb{Z}_{8}, +), M_1 \rangle$.

{\small
{\em Authors' addresses}:
{\em Miroslav Ploščica}, Faculty of Natural Sciences, Šafárik's University, Jesenná 5, 04154 Košice\\
 e-mail: \texttt{miroslav.ploscica@upjs.sk \\ \allowbreak https://ploscica.science.upjs.sk},\\
{\em Radka Schwartzová}, Faculty of Natural Sciences, Šafárik's University, Jesenná 5, 04154 Košice\\
 e-mail: \texttt{radka.schwartzova@student.upjs.sk},\\
{\em Ivana Varga}, Faculty of Natural Sciences, Šafárik's University, Jesenná 5, 04154 Košice\\
 e-mail: \texttt{ivana.varga@student.upjs.sk}.
}

\end{document}